
\documentstyle{amsppt}
\magnification=\magstep1
\NoRunningHeads

\vsize=7.4in

\def\tr{\text{tr }}
\def\det{\text{det }}
\def\dim{\text{dim }}
\def\ubc{\text{ubc }}

\def\Ave{\mathop{\text {Ave}}}
\def\rank{\text{rank }}

\topmatter

\title Unconditional bases and unconditional finite-dimensional
decompositions in Banach spaces
\endtitle
\author
P.G. Casazza
and N.J. Kalton
\endauthor
\address
Department of Mathematics,
University of Missouri,
Columbia, Mo.  65211, U.S.A.
\endaddress
\email pete\@casazza.cs.missouri.edu, mathnjk\@mizzou1.missouri.edu
\endemail
\thanks Both authors were supported  by  NSF
Grant DMS-9201357 \endthanks
\abstract
Let $X$ be a Banach space with an unconditional
finite-dimensional Schauder decomposition
$(E_n)$.  We consider the general problem of characterizing conditions
under which one can construct an unconditional basis for $X$ by forming
an unconditional basis for each $E_n.$   For example, we show that if
 $\sup \dim E_n<\infty$  and $X$ has Gordon-Lewis local unconditional
structure then
$X$
has an
unconditional basis of this type.  We also give an example of a
non-Hilbertian space $X$ with the property that whenever $Y$ is a
closed subspace of $X$ with a UFDD $(E_n)$ such that $\sup\dim
E_n<\infty$ then $Y$ has an unconditional basis, showing that a recent
result of Komorowski and Tomczak-Jaegermann cannot be improved.
\endabstract

\endtopmatter

\document
\baselineskip=14pt

\heading{1. Introduction}\endheading \vskip10pt

Let $X$ be a separable Banach space with an unconditional
finite-dimensional Schauder decomposition (UFDD)  $(E_n).$  It is
well-known that even if for some constant $K$ each $E_n$ has a
$K$-unconditional basis it is still possible that $X$ may fail to
have an unconditional basis.  The first example of this phenomenon was
given in \cite{9} where a twisted sum of two Hilbert spaces $Z_2,$ is
constructed in such a way that it has a UFDD into a two-dimensional
spaces (or a 2-UFDD) $E_n$ but $Z_2$ has no unconditional basis.  Later,
Johnson,
Lindenstrauss and Schechtmann \cite{5} showed that this same example
fails even to have local unconditional structure (l.u.st.).

Recently, Komorowski and Tomczak-Jaegermann \cite{12} proved the
remarkable result that if $X$ has an unconditional basis and is not
hereditarily Hilbertian then it has a subspace $Y$ with
a 2-UFDD and failing local unconditional structure.  This is an
important step in the
resolution
of  the conjecture that a Banach space all of whose subspaces have
local unconditional structure
must be Hilbertian.

Motivated by these results, we investigate here the construction of
unconditional bases or unconditional basic sequences in spaces with a
UFDD.  For convenience let us refer to a UFDD $(E_n)$ as {\it
uniform} if  $\sup_n\dim E_n<\infty$ and as an $N-$UFDD if $\dim
E_n=N$ for all $n.$

Suppose $X$ has an unconditional basis and the property that whenever
$(E_n)$ is a UFDD for $X$ and, for each $n,$ $(f_{ni})_{i=1}^{\dim E_n}$
is a  basis of $E_n$ with unconditional basis constant (ubc) bounded by
some constant $K,$ then $(f_{ni})_{n,i}$ forms an unconditional basis of
$X.$  In Section 2 we prove that this property characterizes the spaces
$\ell_1,\ell_2$ and $c_0.$  A similar property for any UFDD of a
closed subspace characterizes $\ell_2.$

Now suppose $X$ is a Banach space with a uniform UFDD $(E_n).$   Under
these hypotheses we show that (Gordon-Lewis) l.u.st. is equivalent to
the
existence of an unconditional basis for $X$ of the form $(f_{ni})_{n,i}$
where each $(f_{ni})_{i=1}^{\dim E_n}$ is an unconditional basis for
$E_n.$  This provides us with a simple criterion to check whether a
given space with a uniform UFDD has l.u.st.:  compare the earlier
criteria used by Ketonen \cite{10}, Borzyszkowski \cite{2} and
Komorowski \cite{11}.  Using this criterion
we establish a general result on the failure of l.u.st. in twisted
sums.

Finally in Section 4, we give an example to complement the work of
Komorowski and Tomczak-Jaegermann \cite{12}.  We show that there is an
Orlicz sequence space $\ell_F\neq \ell_2$ with the property that
whenever
$(E_n)$
is a uniform UFDD for a closed subspace $X_0$ then one can choose an
unconditional basis $(f_{ni})_{i=1}^{\dim E_n}$ of each $E_n$ so that the
family $(f_{ni})_{n,i}$ is an unconditional basis of $X_0.$  Of course
the space $\ell_F$ is hereditarily Hilbertian; this example shows that
the result of \cite{12} is in a sense best possible.


\heading{2. Preliminary results}\endheading \vskip10pt

Let us say that a UFDD $(E_n)$ is {\it absolute} if there is a constant
$C$ so that if $y_n,x_n\in E_n$ are finitely nonzero and satisfy
$\|y_n\|\le \|x_n\|$ for all $n$ then $\|\sum_{n=1}^{\infty}y_n\|\le
C\|\sum_{n=1}^{\infty}x_n\|.$   We remark that in \cite{1} it is shown
that every FDD of a reflexive subspace of a space with a shrinking
absolute UFDD can be blocked to be a UFDD.
 The following Proposition is trivial:

\proclaim{Proposition 2.1}Suppose $(E_n)$ is an absolute UFDD of a
Banach
space $X$ and that $(f_{ni})_{i=1}^{\dim E_n}$ is an unconditional basis
of $E_n$ so that $\sup_n \ubc (f_{ni})<\infty.$  Then $(f_{ni})_{n,i}$ is
an unconditional basis of $X.$\endproclaim

\proclaim{Proposition 2.2}Let $(E_n)$ be a uniform-UFDD of a Banach
space
$X$ with the property
that whenever we pick an unconditional basis
$(f_{ni})_{i=1}^{\dim E_n}$ of
$E_n$ in such a way that $\sup_n \ubc
(f_{ni})_{i=1}^{\dim E_n}<\infty,$ then  $(f_{ni})_{n,i}$ forms
an unconditional basis of $X.$
Then $(E_n)$ is an absolute UFDD.
\endproclaim

\demo{Proof}Since $\sup\dim E_n<\infty$ we can introduce a
Euclidean norm $\|\,\|_{E_n}$ on each $E_n$ so that for some constant
$K\le \sup_n (\dim E_n)^{1/2}$ we have $\|x\|\le \|x\|_{E_n}\le K\|x\|$
for $x\in E_n.$  We now claim that there is a constant $C$ so that if
$A_n:E_n\to E_n$ is a sequence of operators which are Hermitian
for the Euclidean norm and satisfies
$\sup_n\|A_n\|_{E_n}\le
1$ then for any finitely nonzero sequence $(x_n)$ with $x_n\in E_n,$ we
have
$$ \|\sum_{n=1}^{\infty} A_nx_n\|\le C\|\sum_{n=1}^{\infty}x_n\|.\tag
1$$
Indeed if such an estimate fails, by a simple gliding hump argument we
can construct a single sequence of Hermitian operators $B_n:E_n\to E_n$
with
$\|B_n\|_{E_n}\le 1$ so that (1) fails for $(B_n)$ for any constant
$C.$  But  if we choose $(f_{ni})_{i=1}^{\dim E_n}$ to be an
orthonormal basis for $E_n$ in its Euclidean norm consisiting of
eigenvectors for $B_n$ then, by our hypothesis, the
family
$(f_{ni})_{n,i}$ must
be an unconditional basis for $X$ and this contradicts the failure of (1)
for the sequence $(B_n).$

Now it follows from (1) that if $A_n:E_n\to E_n$ is any bounded sequence
of operators (not necessarily Hermitian) and $x_n\in E_n$ is finitely
nonzero then
$$
\align
 \|\sum_{n=1}^{\infty}A_nx_n\| &\le
2C\sup_n(\|A_n\|_{E_n})\|\sum_{n=1}^{\infty}x_n\|\\
&\le 2CK\sup_n\|A_n\|\|\sum_{n=1}^{\infty}x_n\|
\endalign
$$
where $\|A_n\|$ represents the norm of $A_n$ with respect to the original
norm on $X.$  The proposition now follows immediately.\enddemo

We shall say that an unconditional basis $(e_n)_{n=1}^{\infty}$ for a
Banach space
$X$ has the
{\it shift property} (SP)
if whenever $(x_n)$ is a normalized block basic sequence there is a
constant $C$ so that for any finitely nonzero sequence $(\alpha_n)$ we
have $$ C^{-1}\|\sum_{n=1}^{\infty}\alpha_nx_n\| \le
\|\sum_{n=1}^{\infty}\alpha_nx_{n+1}\| \le
 C\|\sum_{n=1}^{\infty}\alpha_nx_n\|.\tag 2$$

It is easy to see that if $X$ has (SP) then there is a uniform constant
$C$ so that (2) holds for all normalized block basic sequences.  We also
remark that, although our formulation is mildly different, essentially
the same concept was introduced for sequence spaces in \cite{8}.
Precisely the unconditional basis $(e_n)$ has (SP) if and only if the
corresponding sequence space has both the left-shift (LSP) and
right-shift (RSP) properties.  No example of a sequence space with
just one shift property, say (LSP), and not the other is known.

\proclaim{Proposition 2.3}The following properties are
equivalent:\newline
(1) $(e_n)_{n=1}^{\infty}$ has property (SP).\newline
(2) For every blocking $E_n=[e_i]_{i=r_{n-1}+1}^{r_n}$
of
$(e_n)$ and every unconditional basis
$(f_k)_{k=r_{n-1}+1}^{r_n}$ of $E_n$ such that
$\sup_n\ubc (f_k)_{k=r_{n-1}+1}^{r_n}<\infty$ the sequence
$(f_k)_{k=1}^{\infty}$ forms an unconditional basis of $X.$\newline
(3) For every blocking $E_n$  of $(e_n)$ and every
sequence
$(F_n)$ of 2-dimensional spaces so that $F_n\subset E_n$
and every unconditional basis $(f_{2n-1},f_{2n})$ of
$F_n$ with $\sup_n \ubc (f_{2n-1},f_{2n})<\infty$ the sequence
$(f_n)_{n=1}^{\infty}$ is an unconditional basic sequence.
\endproclaim

\demo{Proof}Clearly (1) implies easily that every blocking $(E_n)$ is an
absolute UFDD, and so implies (2) by Proposition 2.1.
(2) implies that every blocking is absolute and so also implies (3) by
Proposition 2.1.

Finally suppose we have (3). Suppose $(x_n)$ is a normalized block basic
sequence. It follows from Proposition 2.2 and (3) that the UFDD
$F_n=[x_{2n-1},x_{2n}]$ is absolute.  Hence, for a suitable constant
$C_0,$ and for any
finitely nonzero sequence $(\alpha_n)$ we have:
$$ C_0^{-1}\|\sum_{n=1}^{\infty}\alpha_{2n-1}x_{2n}\|\le
\|\sum_{n=1}^{\infty}\alpha_{2n-1}x_{2n-1}\|\le
 C_0\|\sum_{n=1}^{\infty}\alpha_{2n-1}x_{2n}\|.$$
In  a similar fashion, considering $G_n=[x_{2n},x_{2n+1}]$ we have a
constant $C_1$ so that:
$$ C_1^{-1}\|\sum_{n=1}^{\infty}\alpha_{2n}x_{2n+1}\|\le
\|\sum_{n=1}^{\infty}\alpha_{2n}x_{2n}\|\le
 C_1\|\sum_{n=1}^{\infty}\alpha_{2n}x_{2n+1}\|.$$
Combining this with the fact that $(x_n)$ is an unconditional basic
sequence shows that $(e_n)$ has the shift property (SP).\enddemo

\proclaim{Theorem 2.4}Let $X$ be a Banach space with an unconditional
basis.  Then the following are equivalent:\newline
(1) $X$ is isomorphic to one of the spaces $\ell_1,\ell_2$ or
$c_0.$\newline
(2) Whenever $(E_n)$ is a UFDD for $X$ and $(f_{ni})_{i=1}^{\dim E_n}$
is an unconditional basis for each $E_n$ with $\sup_n\ubc
(f_{ni})_{i=1}^{\dim E_n}<\infty$ then $(f_{ni})_{n,i}$ is an
unconditional basis for $X.$
\endproclaim

\demo{Proof}$(1)\Rightarrow (2).$
This is obtained by putting together some folklore results.  It follows
easily from the parallelogram law that if $(E_n)$ is a UFDD for
$\ell_2$ then there is a constant $C$ so that if $(x_n)_{n=1}^{\infty}$
is a finitely nonzero sequence with $x_n\in E_n$ then
$$ C^{-1}(\sum_{n=1}^{\infty}\|x_n\|^2)^{1/2}\le
 \|\sum_{n=1}^{\infty}x_n\|\le
 C(\sum_{n=1}^{\infty}\|x_n\|^2)^{1/2}.$$
If $(E_n)$ is a UFDD for $\ell_1$ one obtains the similar inequality
$$ C^{-1}\sum_{n=1}^{\infty}\|x_n\|\le
 \|\sum_{n=1}^{\infty}x_n\|\le
 C\sum_{n=1}^{\infty}\|x_n\|,$$
from the classical argument of Lindenstrauss-Pe\l czy\'nski \cite{13}
that the unconditional basis of $\ell_1$ is unique.  In the case of $c_0$
one obtains
$$ C^{-1}\max_{1\le n<\infty}\|x_n\|\le
 \|\sum_{n=1}^{\infty}x_n\|\le
 C\max_{1\le n<\infty}\|x_n\|.$$
In all cases the UFDD is absolute and Proposition 2.1 gives the
conclusion.

$(2)\Rightarrow (1).$  It follows from Proposition 2.3 that  every
permutation of every unconditional basis has the shift property.  Thus
any unconditional basis $(e_n)$ is a symmetric basis with the (SP) and so
by Proposition
2.3 of
\cite{8} $X$ is isomorphic to one of the spaces $\ell_p$ for $1\le
p<\infty$ or to $c_0.$  Since every unconditional basis is symmetric
this shows that $X$ is isomorphic to one of the three spaces
$\ell_1,\ell_2$ or $c_0$ (cf. \cite{17}).
\enddemo

\proclaim{Corollary 2.5}
Let $X$ be a Banach space with an unconditional
basis.  Then the following are equivalent:\newline
(1) $X$ is isomorphic to  $\ell_2.$ \newline
(2) Whenever $(E_n)$ is a UFDD for a closed subspace $Y$ of $X$ and
$(f_{ni})_{i=1}^{\dim
E_n}$ is an unconditional basis for each $E_n$ with $\sup_n\ubc
(f_{ni})_{i=1}^{\dim E_n}<\infty$ then $(f_{ni})_{n,i}$ is an
unconditional basis for $Y.$ \newline
(3) Whenever $(E_n)$ is a 2-UFDD for a closed subspace $Y$ of $X$ and
$(f_{ni})_{i=1,2}
$ is an unconditional basis for each $E_n$ with $\sup_n\ubc
(f_{ni})_{i=1,2}<\infty$ then $(f_{ni})_{n,i}$ is an
unconditional basis for $Y.$
\endproclaim

\demo{Proof}Clearly (1) implies (2) and (2) implies (3).  For $(3)
\Rightarrow
(1)$ we use Proposition 2.3 to deduce that every unconditional basic
sequence has the shift-property and hence as in Theorem 2.4, $X$ is
isomorphic to $\ell_1$, $\ell_2$ or $c_0.$  Since this property passes to
every closed subspace we obtain $X$ isomorphic to $\ell_2.$
\enddemo

Our final result of the section is that if one can choose an
unconditional basis from a uniform UFDD then it is essentially unique.

\proclaim{Proposition 2.6}Suppose $X$ is a Banach space with a uniform
UFDD $(E_n)$.  Suppose $(f_{ni})_{i=1}^{\dim E_n}$ and
$(g_{ni})_{i=1}^{E_n}$ are normalized unconditional bases for each $E_n,$
so that the whole collections $(f_{ni}),(g_{ni})$ are unconditional bases
of $X$.  Then $(f_{ni})$ and $(g_{ni})$ are permutatively equivalent.
\endproclaim

\demo{Proof}Let $d_n=\dim E_n.$  Let $g_{ni}=\sum_{j=1}^na_{ij}^nf_{nj}.$
It is easy to see that $\inf_n|\det(a_{ij}^n)|>0$ and so for some $c>0$
and each $n,$
there is a permutation $\sigma_n$ of $\{1,2,\ldots,d_n\}$ so that
$|a_{i,\sigma_n(i)}|>c.$

Now by Krivine's theorem, for any finitely nonzero $(\alpha_{ni})$
$$
\|\sum_{n,j}(\sum_i|\alpha_{ni}|^2|a^n_{ij}|^2)^{1/2}f_{nj}\|
\le C\|\sum_{n,i}\alpha_{ni}g_{ni}\|
$$
where $C$ is a suitable constant.
Therefore,
$$ \|\sum_{n,i}\alpha_{ni}f_{n,\sigma_n(i)}\| \le
Cc^{-1}\|\sum_{n,i}\alpha_{ni}g_{ni}\|.$$

Thus the basis $(g_{ni})$ dominates the basis $(f_{n,\sigma_n(i)})$.  By
the same argument, there exist permutations $\tau_n$ of
$\{1,2,\ldots,d_n\}$ so that $(f_{ni})$ dominates $(g_{n,\tau_n(i)}).$
Thus $f_{ni}$ dominates $f_{n,\pi_n(i)}$, where $\pi_n=\sigma_n\tau_n.$.
Iterating $N!$ times where $N=\sup d_n$, since $\pi_n^{N!}$ is always the
identity permutation, this implies that
$(f_{ni})$ and
$(f_{n,\pi(i)})$ are actually equivalent and so $(f_{ni})$ and
$(g_{n\sigma_n(i)})$ are equivalent.
\enddemo

 \heading{3. Spaces with local unconditional structure and a
UFDD}\endheading\vskip10pt

Let $Y$ be a space with an unconditional basis $(e_n)$.  We shall say
that a sequence of finite-dimensional subspaces $(E_n)$ forms a {\it
complemented block UFDD} if there is an increasing sequence of integers
$(p_n)_{n=0}^{\infty}$ with $p_0=0$ so that $E_n\subset
F_n=[e_j]_{j=p_{n-1}+1}^{p_n}$ and a projection $P$ on $Y$ so that
$P(F_n)=E_n.$
If further $\sup\dim F_n<\infty$ then $(E_n)$ is a {\it uniform
complemented block UFDD}.

\proclaim{Lemma 3.1}If $(E_n)$ is a uniform complemented block UFDD then
one can choose an unconditional basis $(f_{nj})_{j=1}^{\dim E_n}$ in each
$E_n$ so that $(f_{nj})_{n,j}$ is an unconditional basis for the closed
linear span $X=\sum_{n=1}^{\infty}E_n.$   Furthermore $(f_{nj})$ is
equivalent in a suitable order to a subsequence of $(e_n).$
\endproclaim

\demo{Proof}  We shall prove the first statement by induction on
$M=\sup_n\dim E_n;$ it is clearly true if $M=1.$ Assume the statement
true whenever
$\sup_n\dim
E_n<M$ and suppose $\sup_n\dim E_n=M.$  We first show that it is possible
to pick normalized vectors
$f_{n1}\in E_n$ so that there is a projection $Q:X\to [f_{n1}]$  with
$Q(E_n)\subset E_n,$ for each $n.$
To see this note that for each $n$ we have
 $$\sum_{k=p_{n-1}+1}^{p_n}\langle Pe_n,e_n^*\rangle =\dim E_n$$
and so there exists $p_{n-1}<k_n\le p_n$ so that
$\alpha_n=\langle Pe_{k_n},e_{k_n}^*\rangle>N^{-1}$ where
$N=\sup(p_n-p_{n-1}).$ Let $f_{n1}=Pe_{k_n}$ and and consider the
projection
$Q:Y\to Y$ defined by $Qy=\sum_{n=1}^{\infty}\alpha_n^{-1}\langle
y,e_{k_n}^*\rangle f_{n1}.$   It is readily verified that $Q$ is bounded.

Let $G_n=E_n\cap Q^{-1}\{0\}.$   Then, after deleting trivial spaces,
$(G_n)$ is a uniform complemented block UFDD with $\sup_n\dim E_n\le
M-1.$   We therefore can pick an unconditional basis
$(f_{nj})_{j=2}^{\dim E_n}$ in each so that $(f_{nj})_{n,j}$ is an
unconditional basis of $X.$

To complete the proof let $H_n=P^{-1}\{0\}\cap F_n.$  Then $(H_n)$ is
also a uniform complemented block UFDD.  It is therefore possible to
extend $(f_{nj})_{j=1}^{\dim E_n}$ to an
unconditional basis $(f_{nj})_{j=1}^{\dim F_n}$ of $(F_n)$ in such a way
that $(f_{nj})_{n,j}$ is an unconditional basis of $Y$.  The final
statement follows from Proposition 2.6.\enddemo

\proclaim{Lemma 3.2}
Let $X$ be a finite-dimensional Banach space.  Suppose $X=E_1\oplus
E_2\oplus\cdots E_n$ with associated projections $Q_j:X\to E_j$
satisfying $$\sup_{|\alpha_j|\le 1}\|\sum_{j=1}^n\alpha_jQ_j\|=K.$$
Suppose $Y$ is a finite-dimensional Banach lattice and that $A:X\to Y$
and $B:Y\to X$ are operators so that  $BA=I_X.$  Then
there
is a finite-dimensional Banach lattice $Z$ with a band decomposition
$Z=Z_1\oplus Z_2\oplus\cdots Z_n$ and operators $A_0:X\to Z$, $B_0:Z\to
X$ with $B_0A_0=I_X,$  $A_0(E_j)\subset Z_j,$  $B_0(Z_j)\subset E_j$ and
$\|A_0\|\|B_0\|\le 2K^2\|A\|^2\|B\|^2.$\endproclaim

\demo{Proof}Consider $Z=Y^n$ with the lattice seminorm
 $$\|(y_1,\ldots,y_n)\|_Z =\sup_{\|y^*\|\le 1}\sum_{j=1}^n\langle
|y_j|,|B^*Q_j^*A^*y^*|\rangle.$$
(Strictly speaking $Z$ should be replaced by its Hausdorff
quotient.)  We define  $A_0:X\to Z$ by $A_0x=(AQ_1x,\ldots,AQ_nx)$ and
$B_0:Z\to X$ by $B_0(y_1,\ldots,y_n)=\sum_{j=1}^n Q_jBy_j$.  Clearly
$B_0,A_0$ satisfy all the required properties except possibly the norm
estimates.

If $x\in X$ and $y^*\in Y^*$ with $\|y^*\|\le 1$ then
$$
\sum_{j=1}^n\langle |AQ_jx|,|B^*Q_j^*A^*y^*|\rangle \le \langle(\sum_{j
=1}^n|AQ_jx|^2)^{1/2},(\sum_{j=1}^n|B^*Q_j^*A^*y^*|^2)^{1/2}\rangle.$$
Now, by Khintchine's inequality,
$$ \|\sum_{j=1}^n|AQ_jx|^2)^{1/2}\|\le
2^{1/2}\Ave_{\epsilon_j=\pm1}\|\sum_{j=1}^n\epsilon_jAQ_jx\| \le
2^{1/2}K\|A\|\|x\|.$$
Similarly,
$$ \|(\sum_{j=1}^n|B^*Q_j^*A^*y^*|^2)^{1/2}\|\le 2^{1/2}K\|A\|\|B\|.$$
It follows that $\|A_0\|\le 2K^2\|A\|^2\|B\|.$

On the other hand if $x^*\in X^*$ with $\|x^*\|\le 1$ and if
$(y_1,y_2,\ldots,y_n)=z\in Z$ then
$$
\align
|\langle B_0z,x^*\rangle| &\le\sum_{j=1}^n |\langle
y_j,B^*Q_j^*x^*\rangle| \\
&\le \sum_{j=1}^m\langle |y_j|,|B^*Q_j^*A^*B^*x^*|\rangle\\
&\le \|B\|\|z\|.
\endalign
$$
Hence $\|A_0\|\|B_0\|\le 2K^2\|A\|^2\|B\|^2.$\enddemo

We now recall that a Banach space $X$ has Gordon-Lewis local
unconditional structure (or l.u.st.)
\cite{4} if
there is a constant $C$ so that whenever $E$ is a finite-dimensional
subspace of $X$ there is a finite-dimensional Banach lattice $Y$ and
operators $A:E\to Y,\ B:Y\to X$ with $\|A\|\|B\|\le C$ and $BA=I_E.$
(A stronger form of local unconditional structure is considered in
\cite{3}).

The following Proposition is established by Johnson, Lindenstrauss
and Schechtman \cite{5}, under the additional assumptions
that $X$ has nontrivial cotype and is complemented in its bidual.

\proclaim{Proposition 3.3}Let $X$ be a Banach space with a UFDD
$(E_n).$  Suppose $X$ has local unconditional structure.  Then
there a Banach space
$Y$ with an
unconditional basis $(e_n)$ so that $Y$ contains $X$ and $(E_n)$ is a
complemented block UFDD.\endproclaim

\demo{Proof}We let $Q_n$ be the natural projection of $X$ onto $E_n$ and
set
$K=\sup_n\|Q_n\|.$  Let
$X_n=\sum_{j=1}^nE_j.$  Using l.u.st. and the preceding
Lemma, there is a constant $C$ so that for each $n$ we can find a
finite-dimensional Banach lattice $Z_n$ with a band decomposition
$Z_n=Z_{n1}\oplus\cdots\oplus Z_{nn}$ and operators $A_n:X_n\to Z_n$,
$B_n:Z_n\to X_n$ so that $B_nA_n=I_{X_n},$ $\|A_n\|\le 1,$ $\|B_n\|\le C$
and
$A_n(E_j)\subset Z_{nj},$  $B_n(Z_{nj})\subset E_j$ for $1\le j\le
n.$

Choose $\epsilon_n>0$ to be a sequence such that
$\sum\epsilon_n<(2C)^{-1}.$  Then by an argument of Johnson \cite{16}
we can find for each $n$ and each $1\le j\le n$ a sublattice $Y_{nj}$ of
$Z_{nj}$ with $\dim Y_{nj}\le d_j$  (independent of $n$) and a map
$L_{nj}:A_n(E_j)\to Y_{nj}$ so that $\|L_{nj}z-z\|\le \epsilon_j\|z\|.$

Now let $Y_n=Y_{n1}\oplus \cdots\oplus Y_{nn}$ and
define $\tilde A_n:X_n\to Y$ by
$\tilde A_nx= \sum_{j=1}^n L_{nj}A_nQ_jx$
Then $\|\tilde A_n-A_n\| \le K\sum_{j=1}^n\epsilon_j,$ so that $\|\tilde
A_n\|\le K.$ Further $\|B_n\tilde A_n-I_{X_n}\| \le \frac12.$
Since $\tilde A_n(E_j)\subset Y_{nj}$ and $B_n(Y_{nj})\subset E_j$ the
operator $D_n=(B_n\tilde A_n)^{-1}$ leaves each $E_j$ invariant. Let
$\tilde B_n=D_nB_n$;  then $\|\tilde B_n\|\le 2C.$

The conclusion from these calculations, after relabelling, is that there
is a constant $C'$ so that for each $n$ there is a Banach lattice $Z_n$
with a band decomposition $Z_{n1}\oplus\cdots \oplus Z_{nn}$ such that
$\dim Z_{nj}\le d_j$ and operators $A_n:X_n\to Z_n$, $B_n:Z_n\to X_n$ so
that $B_nA_n=I_{X_n}$, $\|A_n\|\le 1,$ $\|B_n\|\le C'$ and
$A_n(E_j)\subset Z_{nj}$, $B_n(Z_{nj})\subset E_j$ for $1\le j\le n.$

Let $p_n=\sum_{j=1}^nd_j.$  We can alternatively regard $Z_n$ as
 the space of sequences $(\xi_j)$ so that $\xi_j=0$ for $j>p_n,$ with an
associated norm.  We can further suppose that the canonical basis vectors
$(e_j)_{j=1}^{p_n}$ are normalized and that
$Z_{nj}=[e_k]_{k=p_{j-1}+1}^{p_j}.$

Let $\Cal U$ be a nonprincipal ultrafilter on the natural numbers $\bold
N$.  Define a norm $\|\,\|_Z$ on the space $c_{00}$ of all finitely
nonzero sequences by $\|\xi\|_Z=\lim_{\Cal U}\|\xi\|_{Z_n}.$
Let $Z$ be the completion of $c_{00}$ for this norm.

Let $X_0$ be the linear span of all $(E_n)$ in $X$.  We can define an
operator $A:X_0\to c_{00}$ by $Ax=\lim_{\Cal U}A_nx$ and similarly
$B:c_{00}\to X_0$ by $Bx=\lim_{\Cal U}B_n\xi.$
It is clear that $\|A\|\le 1$ and $\|B\|\le C'.$  It is easy to verify
that $A$ isomorphically embeds $X$ into $Z$ in such a way that $A(E_n)$
is a complemented block UFDD.\enddemo

\demo{Remark}In \cite {5} it is further claimed in Remark 2 that (under
their
additional hypotheses) if $(E_n)$ is a {\it uniform} UFDD then we can
choose
$Y$ so that $(E_n)$ is a {\it uniform} complemented block UFDD.  This of
course would imply by Lemma 3.1, that one could find an unconditional
basis for
$X$
by picking a basis of each $E_n.$  However, no proof of Remark 2 is given
and the natural proof does not appear to work.  We shall see, however,
that the claim of Remark 2 in \cite{5} is nonetheless correct, but the
proof is rather circuitous.  \enddemo

Let us now fix $H$ as a complex $N$-dimensional Hilbert space, where
$N\ge
2.$   If $A\in\Cal L(H)$ we denote its trace by $\tr A$ and its spectral
radius by $r(A).$
 We say
that a subalgebra $\Cal A$ of $\Cal L(H)$ is {\it triangular} if every
$A\in\Cal A$ is of the form $A=\lambda I +S$ where $S$ is nilpotent.
This is equivalent to requiring that $r(A-\frac1N(\tr A)I)=0$ for all
$A\in \Cal A.$ The folloing elementary lemma is very well-known and we
include its proof only for reference.

\proclaim{Lemma 3.4}
If $\Cal A$ is a triangular subalgebra of $\Cal L(H)$ then
 there is an orthonormal basis $(e_j)_{j=1}^N$ so
that every
$A\in\Cal A$ with $\tr A=0$ is upper triangular i.e.
$(e_j,Ae_k)=0$ whenever $j\le k.$ \endproclaim

\demo{Proof}The subset $\Cal A_0=\{A\in\Cal A:\tr A=0\}=\{A\in\Cal A:\det
A=0\}$ is an ideal of $\Cal A.$  It suffices to construct an increasing
sequence of subspaces $(E_k)_{0\le k\le N}$ with $\dim E_k=k$ so that
$\Cal A_0(E_k)\subset A_{k-1}$ for $1\le k\le N.$  Then we can construct
an
orthonormal basis $(e_k)_{k=1}^N$ with $e_k\in E_{N-k+1}$ for $1\le k\le
N.$  Suppose then $E_0=\{0\}$ to start the induction.  Now suppose $1\le
k\le N$ and $E_{k-1}$ has been constructed.  Choose $x\notin E_{k-1}$ to
minimize the dimension of $(\Cal A_0x +E_{k-1})/E_{k-1}.$  If this
dimension is zero then $\Cal A_0x\subset E_{k-1}$ and we let
$E_k=[x,E_{k-1}].$  Otherwise there exists $S\in\Cal A_0$ so that
$Sx\notin E_{k-1}.$   But then $\Cal A_0Sx+E_{k-1}$ contains $Sx$ by
minimality.  Hence there exists $T\in\Cal A_0$ with $Sx-TSx\in E_{k-1}.$
However $(I-T)$ is invertible in $\Cal A$ and $E_{k-1}$ is $\Cal
A$-invariant so that $Sx\in E_{k-1}$ and we have a
contradiction.\enddemo

Now suppose that $\Cal C$ is a compact subset of  $\Cal
L(H)$ which contains the identity $I=I_H.$ Let $\|\Cal C\|
=\sup\{\|A\|:A\in \Cal C\}\ge 1.$  We define for
$m\in\bold N,$
the set $\Cal C^{(m)}$ to be the set of all operators $T\in\Cal L(H)$ of
the form
$$T=\sum_{j_1=1}^m\sum_{j_2=1}^m\cdots\sum_{j_m=1}^m\alpha_{j_1,\ldots,j_m}
A_{j_1}\ldots A_{j_m}$$
where $A_1,A_2,\ldots,A_m\in\Cal C$ and $|\alpha_{j_1,\ldots,j_m}|\le
2^m.$
Since $I \in \Cal C$ the sets $\Cal C^{(m)}$ are increasing compact sets
and $\cup_m\Cal C^{(m)}$ is the algebra generated by $\Cal C.$
  If $T\in\Cal C^{(m)}$
then
$\|T\|\le (2m\|\Cal C\|)^m.$

\proclaim{Lemma 3.5}Suppose $\delta>0,M\ge 1,m,N\in\bold N$ with $N\ge
2.$ Then
there exists $p\in\bold N$ so that $p=p(M,N,m,\delta)$ has the
following property:  Suppose $\Cal C$ satisfies the above
conditions with
$\|\Cal C\|\le M.$  Suppose
$S\in\Cal C^{(m)}$ and
$r(S-\frac1N(\tr
S) I)=\delta>0.$  Then there exists a projection $P$ with
$0<\rank P<N$ and $P\in\Cal C^{(p)}.$
\endproclaim

\demo{Proof}
Let $(\lambda_j)_{j=1}^N$ be the (complex) eigenvalues of $S$ repeated
according
to algebraic multiplicity.  We have $\max_{j,k}|\lambda_j-\lambda_k|\ge
\delta.$ It follows that we can reorder them so that for some $s$ with
$1\le s\le N-1$ we have $|\lambda_j-\lambda_k|\ge \frac{\delta}{2N}$
whenever $1\le j\le s$ and $s+1\le k\le N.$
Let $T=\prod_{j=1}^k(S-\lambda_jI).$  For each $j$ we have
$|\lambda_j|\le (2mM)^m$ and on multplying out one has $T\in \Cal
C^{(q)}$ where $q$ depends only on $m,N$ and $M.$
Let
$\mu_k=\prod_{j=1}^s(\lambda_k-\lambda_j)$ so that $\mu_k=0$ if $1\le
k\le s$ and $2^N(2mM)^{mN}\ge |\mu_k|\ge (\delta/(2N))^N$ if $s+1\le
k\le N.$ Next let
$W=\prod_{k=s+1}^N(T-\mu_kI).$  Let $\gamma=\prod_{k=s+1}^N(-\mu_k)$. It
is easily seen that
$P=\gamma^{-1}W$ is a projection and that $0<\rank P<N$.  From the
obvious upper bound on $\gamma^{-1}$, we
obtain immediately that $P\in\Cal C^{(p)}$ where $p$ depends
only on
$m,N,M,\delta.$
\enddemo

 The next estimate is crude and
can doubtless be improved.

\proclaim{Lemma 3.6}Suppose $H$ is a complex $N$-dimensional Hilbert
space and that $\Cal A$ is a triangular subalgebra of $\Cal L(H).$
 Let $(A_k)$ be a sequence in $\Cal L(H)$ with each $A_k$
non-invertible such that $\sum_{n=1}^{\infty}A_k$ converges
unconditionally to $I=I_H.$  Let
$$ \sup_{|\alpha_k|\le 1}\|\sum_{k=1}^{\infty}\alpha_kA_k\|=M.$$
Then
$$ \sup_{|\alpha_k|\le 1}d(\sum_{k=1}^{\infty}\alpha_k
A_k,\Cal A)\ge 2^{-3N^2}(N!)^{-N}M^{1-N^2}.$$
\endproclaim

\demo{Proof}
Let $b= \sup_{|\alpha_k|\le 1}d(\sum_{k=1}^{\infty}\alpha_k
A_k,\Cal A).$
First observe that since
$\Cal A$ is triangular
 we
can choose an orthonormal basis $(e_j)_{j=1}^N$ so that, when represented
as matrices,  each
$B\in\Cal A$
is of the form $B=\lambda I +S$ where $S$ has an upper triangular matrix
i.e. $S=(s_{ij})$ where $s_{ij}=0$ if $i\le j.$

Next, note that there exists $B\in\Cal A$ so that $\|I-B\|\le b.$
If $B=\lambda I +S$ then $\tr B=N\lambda$ and so $|\lambda-1|
\le b.$    Then $\|I-S\| \le (b+\tau),$ where $\tau=|\lambda|.$
Clearly $\|S\| \le 1+2b+\tau\le 4M.$  Hence on expanding since $S^N=0$
we obtain $1 \le 2^N(b+\tau)(4M)^{N-1}=
2^{3N-2}(b+\tau)M^{N-1}.$

We now estimate $\tau.$  Let $A_n=\{a^n_{ij}\}_{i,j=1}^N.$  Then
$$
\align
\frac1N|\tr A_n|&\le \frac1N\sum_{i=1}^N|a^n_{ii}| \\
                &\le
(\prod_{i=1}^N|a^n_{ii}|)^{1/N}+\max_{i>j}|a^n_{ii}-a^n_{jj}|\\
&\le
(\prod_{i=1}^N|a^n_{ii}|)^{1/N}+\sum_{i>j}|a^n_{ii}-a^n_{jj}|.
\endalign
$$
Notice that for fixed $i,j$ we have
$\sum_{n=1}^{\infty}|a^n_{ii}-a^n_{jj}| \le 2b.$  Hence on summing we
have
$$ \tau \le \sum_{n=1}^{\infty}|\prod_{i=1}^Na^n_{ii}|^{1/N} + N^2b.$$

Now since $\det A_n=0$ we have that
$$ \prod_{i=1}^N|a^n_{ii}| \le \sum_{\sigma\in
\Pi'}\prod_{i=1}^N|a^n_{i,\sigma(i)}|$$
where $\Pi'$ is the collection of all permutations other than the
identity of
$\{1,2,\ldots,N\}.$
Hence
$$ \sum_{n=1}^{\infty}
 (\prod_{i=1}^N|a^n_{ii}|)^{1/N} \le
\sum_{\sigma\in\Pi'}\sum_{n=1}^{\infty}(\prod_{i=1}^N|
a^n_{i,\sigma(i)}|)^{
1/N}.$$

Let us fix $\sigma\in\Pi'.$  Then
$$\sum_{n=1}^{\infty}(\prod_{i=1}^N|a^n_{i,\sigma(i)}|)^{
1/N} \le \prod_{i=1}^N(\sum_{n=1}^{\infty}|a^n_{i,\sigma(i)}|)^{1/N}.$$

Now if $i>\sigma(i)$ then
$$ \sum_{n=1}^{\infty}|a^n_{i,\sigma(i)}| \le b$$
so that since $\sigma$ is not the identity, we obtain an upper estimate
$$\sum_{n=1}^{\infty}(\prod_{i=1}^N|a^n_{i,\sigma(i)}|)^{1/N} \le
M^{1-1/N}b^{1/N}.$$

Summing over all such permutations and combined with our previous
estimates we finally obtain:
$$ \tau \le N^2b + (N!)M^{1-1/N}b^{1/N}.$$
It follows that
$$ b+\tau \le 4(N!)M^{1-1/N}b^{1/N}$$
and the lemma follows.\enddemo

\proclaim{Lemma 3.7}Suppose $2\le N\in\bold N$ and $M\ge 1$  Then
there is an integer $p=p(N,M)$ with the following property:
suppose
$H$ is a complex
$N$-dimensional Hilbert space and
let $(A_k)$ be a sequence in $\Cal L(H)$ with each $A_k$
non-invertible such that $\sum_{k=1}^{\infty}A_k$ converges
unconditionally to $I=I_H.$  Let $\Cal C$ be the collection of all
operators of the form $\sum_{k=1}^{\infty}\alpha_kA_k$ where
$|\alpha_k|\le 1$ and suppose $\|\Cal C\|\le M.$  Then there is a
projection $P$ with $0<\rank P<N$ and $P\in\Cal C^{(p)}.$
\endproclaim

\demo{Proof}We argue by contradiction.  Suppose for some $M,N$ the result
is false.  Then we can find a sequence of such expansions
$I=\sum_{k=1}^{\infty}A_{nk}$ so that the associated compact sets
$\Cal C_n=\{\sum_{k=1}^{\infty}\alpha_kA_{nk}:|\alpha_k|\le 1\}$
satisfy $\|\Cal C_n\|\le M$ and that there is no nontrivial projection in
$\Cal C_n^{(n)}.$   By passing to a subsequence we can further suppose
that $\Cal C_n$ converges in the Hausdorff metric to a compact set $\Cal
C.$  It  then follows from Lemma 3.5  that for each $p$ we must have
$$ \lim_{n\to\infty}\sup_{S\in\Cal C_n^{(p)}}r(S-\frac{\tr S}N I)=0.$$
Since all these quantities are continuous it follows that  if $S\in\Cal
C^{(p)}$ for any $p,$ we have $$r(S-\frac1N (\tr S)I)=0$$  and
so the
algebra $\Cal A$ generated by $\Cal C$ is triangular.  But the preceding
Lemma 3.6 now implies that
$$ \inf_n \sup_{A\in\Cal C_n}d(A,\Cal A)>0$$
which contradicts the fact that $\Cal C_n$ converges in the Hausdorff
metric to $\Cal C\subset\Cal A.$\enddemo

\proclaim{Theorem 3.8}Let $X$ be a real or complex Banach space with
local unconditional structure.  Suppose that $X$ has a uniform UFDD
$(E_n)$. Then there is an
unconditional basis $(f_{nj})_{j=1}^{\dim E_n}$ of each $E_n$ so that
$(f_{nj})_{n,j}$ is an unconditional basis for $X.$\endproclaim

\demo{Proof}We first prove the complex case.  We shall prove the formally
weaker statement that if $X$ has l.u.st. and an $N$-UFDD $(E_n)$ then
there
is a bounded projection $Q$ on $X$ so that $Q(E_n)\subset E_n$ for each
$n$ and $0<\dim Q(E_n)<N$ for each $N.$  Once this is proved the result
follows simply by induction on $\sup_n\dim E_n.$

We first note that it is possible by Proposition 3.3 to regard $(E_n)$ as
a complemented block UFDD in a Banach space $Y$ with unconditional basis
$(e_n).$  We suppose that $E_n\subset [e_k]_{k=r_{n-1}+1}^{r_n}$ where
$r_0<r_1<\cdots.$   Let
$P:Y\to X$ be the associated projection. Let
$H$ be an
$N$-dimensional Hilbert space and suppose for each $n,$ $V_n:H\to E_n$ is
an isomorphism satisfying $\|V_n\|\|V_n^{-1}\|\le \sqrt N.$

Letting $(e_k^*)$ be the biorthogonal functions for the basis, we define
for $r_{n-1}+1\le k\le r_n,$  $A_{nk}:H\to H$ by
$A_{nk}=V_n^{-1}P(e_k^*\otimes e_k) V_n.$  Let $\Cal
C_n=\{\sum_{k=r_{n-1}=1}^{r_n}\alpha_kA_k:\ |\alpha_k|\le 1\}.$
Then note that $\sup_n\|\Cal C_n\|<\infty$ and so there is an integer $p$
such each $\Cal C_n^{(p)}$ contains a projection $R_n,$ with $0<\rank
R_n<N.$

Next observe that if $(B_n)_{n\in\bold N}$ is any sequence in $\Cal
L(H)$ with $B_n\in \Cal C_n$ then the operator $B$ defined on $X$ by
$Bx=V_nB_nV_n^{-1}x$ for $x\in E_n$ is bounded.  It follows easily that
the same statement is true if $B_n\in\Cal C_n^{(q)}$ for fixed $q.$
Hence the operator $Q:X\to X$ defined by $Qx=V_nR_nV_n^{-1}$ for $x\in
E_n$ is bounded and the proof is complete in the complex case.

We now turn to the real case.  Let $Q_n:X\to E_n$ be the natural
projections.  We complexify $X$ to a space $\tilde X,$ by norming
$(x+iy)$ for
$x,y\in X$ by
$$ \|x+iy\|_c = \sup_{0\le \theta_n\le
2\pi}\|\sum_{k=1}^{\infty}(Q_kx\cos\theta_k +Q_ky\sin\theta_k)\|.$$
Now the subspaces $\tilde E_n=E_n+iE_n$ form a UFDD for $\tilde X$ and so
we can pick an unconditional basis $(\phi_{nj})_{j=1}^{\dim E_n}$ in each
$E_n$ so that $(\phi_{nj})_{n,j}$ is an unconditional basis of $\tilde
X.$   Next let $Y$ be the underlying real space for $\tilde X=X\oplus X.$
Then $Y$ has an unconditional basis $(\phi_{nj},i\phi_{nj})_{n,j}.$  Now
the original $(E_n)$ is a uniform complemented block UFDD in $Y$ with
this basis and so we complete the proof by applying Lemma 3.1.\enddemo

Let us give a sample application of this result.  Let $\omega$ be the
space of all sequences.
  Suppose
$X$ is a super-reflexive (K\"othe) sequence space (so that the canonical
basis vectors $(e_n)$ form a 1-unconditional basis of $X$) and let
$\Omega:X\to\omega$ be a  (homogeneous) centralizer i.e. a map
satisfying, for a suitable constant $\Delta$:\newline
(1) $ \Omega(\alpha x) =\alpha \Omega(x)$ for $\alpha\in\bold R$ and
$x\in X.$\newline
(2) $\|\Omega(ux)-u\Omega(x)\|_X \le
\Delta\|u\|_{\infty}\|x\|_X$ for $x\in X$ and
$u\in\ell_{\infty}$.\newline
See \cite{6} and \cite {7} for discussion and examples.  The
simplest examples are those discussed in \cite{9} of maps
$$ \Omega(x)(n) = x(n) f(\frac{\log |x(n)|}{\|x\|_X})$$
where $f:\bold R\to\bold R$ is a Lipschitz function (here we interpret
the right-hand side as $0$ if $x(n)=0.$)  In the case
$f(t)=t,$ and $X=\ell_2$ one recovers the space $Z_2$ studied in
\cite{9} and \cite{5}.

We can now form the twisted sum $Y=X\oplus_{\Omega}X$ of all pairs
$(x,y)$ in $\omega\times \omega$ such that
$$ \|(x,y)\|_{\Omega}= \|x\|_X +\|y-\Omega(x)\|_X<\infty.$$
This is a quasinorm, but is equivalent to a norm, since the space $X$ is
superreflexive (as in \cite{9}). The space $Y$ is then a reflexive
Banach space with a
2-UFDD
$(E_n)$
where $E_n$ is the span of $(e_n,0)$ and $(0,e_n).$
The vectors $(0,e_n)$ span a closed subspace $X_0$ of $Y$ isomorphic to
$X$ and the quotient space $Y/X_0$ is also isomorphic to $X$.   It
may be shown that $X_0$ is complemented in $Y$, so that $Y$ splits as a
direct sum $X\oplus X$, if and only
if there is a linear centralizer $L:X\to \omega$ (i.e. $Lx=bx$ for some
$b\in\omega$) so that $\|Lx-\Omega x\|_X \le C\|x\|_X$ for all $x\in
X.$  Such
twisted sums arise very naturally as derivatives of complex interpolation
scales of sequence spaces.  If $Z_0,Z_1$ are are two
super-reflexive
sequence spaces and $Z_{\theta}=[Z_0,Z_1]_{\theta}$ for $0<\theta<1$ is
the usual interpolation space by the Calder\'on method, one can define a
derivative $dX_{\theta}$ which is a twisted sum
$X_{\theta}\oplus_{\Omega}X_{\theta}$ which splits if and only if
$Z_1=wZ_0$ for some weight sequence $w=(w(n))$ where $w(n)>0$ for all
$n.$   These remarks follow easily from the methods of
\cite{6}.

Our main conclusion here is that twisted sums of this type have
l.u.st. if
and only if they split as a direct sum.  This extends the special case of
$Z_2$ given in \cite{5}.

\proclaim{Theorem 3.9}If $X$ is a super-reflexive sequence
space
and $\Omega:X\to\omega$ is a centralizer on $X$.  Let
$Y=X\oplus_{\Omega}X.$ Then the following
are equivalent:\newline (1) $Y$ is isomorphic to $X\oplus X.$\newline
(2) $Y$ has local unconditional structure.\newline
(3) $Y$ has an unconditional basis.\newline
(4) The subspace $X_0$ is complemented in $Y.$\newline
(5) There exists $b\in\omega$ and $C>0$ so that $\|\Omega(x)-bx\|_X \le
C\|x\|_X$ for all $x\in X.$\endproclaim

 \demo{Proof}We have already observed the equivalence of (4) and (5).
Clearly $(4) \Rightarrow (1)\Rightarrow (3) \Rightarrow (2).$  It
remains only to show that $(2)\Rightarrow (4).$

Let us first remark that we can assume that the canonical basis $(e_n)$
of $X$ is normalized; we can also assume that $X$ is $p$-convex and
$q$-concave with constant one, for suitable $p>1$ and
$\frac1p+\frac1q=1.$ We note that
$Y$ is super-reflexive and has a
2-UFDD
$(E_n)$ with the
property that for suitable $x_n\in E_n$ the unconditional basic sequence
$(x_n)$ is equivalent to the canonical basis of the sequence space $X$
and induced unconditional basis $(y_n)$ of the quotient $Y/[x_n]$ is also
equivalent to the canonical basis of $X.$

If $Y$ has l.u.st. then by Theorem 3.8,  we can pick a normalized
basis of
$E_n$ say $(f_n,g_n)$ so that $(f_n,g_n)_{n=1}^{\infty}$ is an
unconditional basis of $Y.$  We may suppose that $x_n=a_n f_n+b_ng_n$
where $|a_n|\ge b_n\ge 0.$  If we consider the dual basis $(f_n^*,g_n^*)$
then the sequence $x_n^*=b_n f_n^*-a_ng_n^*$ must be equivalent to the
canonical basis of the dual sequence space $X^*$.

It will now be convenient to switch to sequence space language.  Let $W$
be a sequence space so that if $\xi\in W$ then $\|\xi\|_W$ is equivalent
to $\|\sum_{k=1}^{\infty}\xi(k)f_k\|$ and similarly let $Z$ be a sequence
space so that $\|\xi\|_Z$ is equivalent to
$\|\sum_{k=1}^{\infty}\xi(k)g_k\|$.  Since $Y$ is super-reflexive we can
assume that both $W$ and $Z$ are $p$-convex and $q$-concave with constant
one (possibly changing the original choice of $p,q$), and that the
canonical bases are normalized in both $W$ and $Z.$

It is now easy to see that for a suitable constant $C$ we have the
inequalities
$$ \frac1C\|\xi\|_X \le \max(\|\xi\|_W,\|b\xi\|_Z)\le C\|\xi\|_X$$
and
$$ \frac1C\|\xi^*\|_{X^*}\le \max(\|b\xi^*\|_{W^*},\|\xi^*\|_{Z^*})\le
C\|\xi^*\|_{X^*}$$
whenever $\xi,\xi^*\in c_{00}.$  Note also that $b\in\ell_{\infty}$.  We
will show that these inequalities imply that $W$ and $Z$ both coincide up
to equivalence of norm with $X$ and hence that the basic sequences
$(f_n),(g_n)$ and $(x_n)$ are all equivalent.  This suffices to show that
$[x_n]$ is indeed complemented in $Y,$ i.e. $X_0$ is complemented in $Y.$

Before proceeding we will need a lemma:
\enddemo

\proclaim{Lemma 3.10}Suppose $V$ is a $p$-convex sequence
space with $1<p<\infty,$ and that $0\le \xi,\xi^*\in c_{00}$ with
$\langle
\xi,\xi^*\rangle
=\|\xi\|_V=\|\xi^*\|_{V^*}=1.$  Suppose further that $0\le \eta\in
c_{00}$ with
$\langle \eta,\xi^*\rangle \ge1$ and $\|\eta\|_{V}=M.$  Then, if
$\frac1p+\frac1q=1,$
$\langle\min(\xi,\eta),\xi^*\rangle \ge q^{-1}M^{-q}.$\endproclaim

\demo{Proof of Lemma}Note that if $t\ge 0,$
$$ \langle \max(\xi,t\eta),\xi^*\rangle \le (1+M^pt^p)^{1/p}$$
and so
$$ \langle \min(\xi,t\eta),\xi^*\rangle \ge 1+t-(1+M^pt^p)^{1/p}.$$
Now let $t=M^{-q}\le 1$; the lemma follows by elementary
estimates.\enddemo

\demo{Proof of the Theorem: (2) implies (4)}
We observe first that if $\xi\in c_{00}$ then $\|\xi\|_W\le C\|\xi\|_X
\le C^2\|\xi\|_Z.$

Suppose first that $M$ is chosen so large that
 $ 2qC^{4q+2} + M(1-\frac1{2^q})^{1/q}<M.$
Let
$\beta>0$ be chosen so that $\beta<\min(C^{-2},(M+1)^{-1}).$
 We split
$\bold N$ as
$\Cal M'\cup\Cal M$ where $\Cal M=\{n:b_n\le \beta\}$ and $\Cal
M'=\{n:b_n>\beta\}.$  First suppose $\xi\in c_{00}$ is supported on
$\Cal M'.$  Then
$$ \max(\|\xi\|_W,\|\xi\|_Z)\le C\beta^{-1}\|\xi\|_X$$
and if $\xi^*\in c_{00}(\Cal M'),$ then
$$ \max(\|\xi^*\|_{W^*},\|\xi^*\|_{Z^*})\le C\beta^{-1}\|\xi^*\|_{X^*}.$$
These inequalities show that on $c_{00}(\Cal M')$ the spaces $X,W,Z$
coincide.

Now let $\kappa_n$ be the supremum of $\|\xi\|_Z$ subject to $\xi\in
c_{00}(\Cal M),$  $\|\xi\|_W= 1$ and $\xi$ has support of cardinality
at most $n.$  Then $\kappa_1=1$ and $\kappa_{n+1}\le \kappa_n +1.$
We will show by induction that $\kappa_n\le M$ for all $n.$

Suppose $\kappa_{n-1}\le M,$ and $\kappa_n>M.$  Then there exists $\xi\ge
0$ with support of exactly $n$ so that $\xi\in c_{00}(\Cal M)$,
$\|\xi\|_W=1$ and $\|\xi\|_Z>M.$   Now $\|\xi\|_Z\le M+1$ so that
$\|b\xi\|_Z \le \beta(M+1)<1.$  Hence we must have $\|\xi\|_X\le
C.$   Pick $0\le \xi^*$ with the same support so that $\langle
\xi,\xi^*\rangle =1$ and $\|\xi^*\|_{W^*}=1.$   Then $\|\xi^*\|_{X^*}\ge
C^{-1},$ but
$\|b\xi^*\|_{W^*}\le \beta<C^{-2},$ so that $\|\xi^*\|_{Z^*}\ge C^{-2}.$
It follows that we can pick $0\le \eta$ again with the same support,
so that $\|\eta\|_Z \le C^2$ and $\langle\eta,\xi^*\rangle =1.$  We then
have $\|\eta\|_W\le C^4,$ and we can apply the lemma to see that
$$ \langle \min(\xi,\eta),\xi^*\rangle \ge q^{-1}C^{-4q}.$$
Now let $\zeta(n)= \xi(n)$ whenever $\xi(n) \le 2qC^{4q}\eta(n)$ and
$\zeta(n)=0$ otherwise.  It follows easily that
$\langle\zeta,\xi^*\rangle\ge \frac12$ (so that $\|\zeta\|_W\ge
\frac12$) and
$\|\zeta\|_Z
\le
2qC^{4q}\|\eta\|_X \le 2qC^{4q+2}.$

Now
$$
\align
\|\xi\|_Z &\le \|\zeta\|_Z +\kappa_{n-1}\|\xi-\zeta\|_W \\
 &\le 2qC^{4q+2} + M(1-\|\zeta\|_W^q)^{1/q}\\
 &\le 2qC^{4q+2} + M(1-\frac1{2^q})^{1/q}\\
 &\le M.
\endalign
$$
This contradiction yields the result that $\kappa_n\le M$ for all $n$ and
hence the theorem.\enddemo

\demo{Remark}The most natural case of Theorem 3.9 is when
$X=\ell_2$ so that $Y$ is a ``twisted Hilbert space'', i.e. $Y$ has a
Hilbertian subspace $X_0$ so that $Y/X_0$ is HIlbertian.  The result
suggests the conjecture that every twisted Hilbert space with an
unconditional basis is a Hilbert space.\enddemo

\heading{4. An example}\endheading\vskip10pt

In this final section we construct an explicit example of a
non-Hilbertian
Orlicz sequence space where every closed subspace with a uniform UFDD
has local unconditional structure.
In \cite{12} Komorowski and Tomczak-Jaegermann show that if a
Banach space with an unconditional basis is not hereditarily Hilbertian
then it has a closed subspace failing l.u.st. but with a uniform UFDD.

We define a function $G$ on $[0,\infty)$ by $G(0)=0$ and:
$$
G(x)=
\cases
x(1-\frac12\log x) & 0<x\le 1\\
 \sqrt x & x>1.
\endcases
$$
Note that $G$ is differentiable on $(0,\infty)$ and
$$
G'(x)=
\cases
\frac12(1-\log x) & 0<x\le 1\\
\frac1{2\sqrt x} & x>1.
\endcases
$$

\proclaim{Lemma 4.1}For any $1\le p<\infty,$
 whenever $(a_n)_{n=0}^{\infty}$
is a sequence with $0\le a_n\le 1$ and $(\sum_{n=1}^{\infty}a_n^p)^{1/p}
\le 1$ then for any sequence $(t_n)$ with $t_n\ge 0$ we have:
$$\sum_{n=1}^{\infty}G(a_nt_n) \le G(\sum_{n=1}^{\infty}a_nt_n) +p
\sum_{n=1}^{\infty}t_nG(a_n).$$
\endproclaim

\demo{Proof}
It clearly suffices to prove the inequality for a finite sequence
$(a_1,\ldots,a_N)$ of strictly positive numbers.

Suppose $\lambda>p$.
Consider the function $$\Phi(t_1,\ldots,t_n) =
\sum_{n=1}^NG(a_nt_n)-
G(\sum_{n=1}^Na_nt_n)-\lambda\sum_{n=1}^Nt_nG(a_n),$$
defined on the positive cone $\{t:t_i\ge 0,\ 1\le i\le N\}.$
Note first that if $a_nt_n\ge 1$ then $G(a_nt_n)=a_n^{1/2}t_n^{1/2}\le
a_nt_n\le t_nG(a_n).$  It follows quickly that $\Phi$ is bounded above by
its maximum on the set $\{t:a_it_i\le 1,\ 1\le i\le N\}.$

Let $\Phi$ attain its maximum at the point $s$ where $0\le s_i\le
a_i^{-1}$ for $1\le i\le N.$  Let $S=\sum_{j=1}^Na_js_j.$ For any index
$j$ such that
$s_j>0$ we have
$$ a_jG'(a_js_j) - a_jG'(S) -\lambda G(a_j)=0.$$
Since $a_js_j\le 1$ this simplifies to
$$ \frac12(1-\log a_js_j)-G'(S) =\lambda (1-\frac12\log a_j) $$
and hence to
$$ \log a_js_j+2G'(S)= 1-2\lambda +\lambda\log a_j.$$

Assume that the set $J$ of indices such that $s_j>0$ is nonempty.  Then
taking exponentials and summing
$$ e^{2G'(S)}S = e^{1-2\lambda}\sum_{j\in J}a_j^{\lambda}.$$
Now if $S\ge 1$ then $e^{2G'(S)}S\ge S.$  If $S<1$ then $e^{2G'(S)}S=
e.$  In either case we deduce that
$$ \sum_{j\in J}a_j^{\lambda} \ge e^{2\lambda-1}>1.$$
and this contradicts the conditions on $(a_1,\ldots,a_n).$
Now since $J$ is empty the maximum is attained at the origin and is $0.$
Since $\lambda>p$ is arbitrary the lemma is proved.\enddemo

Now let $F$ be the Orlicz function defined by $F(0)=0$ and
$$
F(x)=
\cases
x^2(1-\log x) & 0<x\le 1\\
x             & x>1.
\endcases
$$
The function $F$ is convex for $x\le e^{-1/2}$ so that $F$ is equivalent
at $0$ to a convex function.  We will consider the Orlicz sequence space
$\ell_F.$  The norm defined in the usual way $$\|x\|_F
=\sup\{t>0:\sum_{n=1}^{\infty}F(x(n)/t)\le 1\}$$ is, strictly
speaking only a quasi-norm but is equivalent to a norm.
Note that
$\ell_F$ is reflexive and has cotype
2 and type
$p$ for any
$p<2;$ these facts are easily computed from the function $F.$  Clearly
$\ell_F\subset
\ell_2$ and
$\|x\|_2\le
\|x\|_F$ for all $x\in\ell_F.$
We will also consider the modular $\Lambda$ defined on $\ell_F$ by
$$\Lambda(x) =\sum_{n=1}^{\infty}F(|x(n)|).$$

It will be convenient to introduce for $0\le a\le 1$ the Orlicz functions
$F_a$ where $F_a(0)=0$ and
$$F_a(x)=
\cases
x^2(1-a\log x) & 0<x\le 1\\
x & x>1.
\endcases
$$
If $(a_n)$ is any sequence with $0\le a_n\le 1$ we will consider the
Orlicz modular space  (or Orlicz-Musielak space) $\ell(F_{a_n})$ of all
sequences
$x(n)$ so that $\sum_{n=1}^{\infty}F_{a_n}(|x(n)|)<\infty$ again with the
(quasi-)norm defined in the usual way.

\proclaim{Lemma 4.2}If $M\ge 1$ there is a constant $K=K(M)$ with the
following property.  Whenever $(x_n)_{n\in \Cal M}$ and $(y_n)_{n\in
\Cal M}$ are two (finite or infinite)
$M-$unconditional basic sequences satisfying the conditions:
$$ \sup_{n\in
\Cal M}\max(\frac{\|x_n\|_F}{\|y_n\|_F},\frac{\|y_n\|_F}{\|x_n\|_F})\le
M$$
$$ \sup_{n\in
\Cal M}\max(\frac{\|x_n\|_2}{\|y_n\|_2},\frac{\|y_n\|_2}{\|x_n\|_2})\le
M$$ then $(x_n)_{n\in \Cal M}$ and $(y_n)_{n\in \Cal M}$ are
$K$-equivalent.
\endproclaim

\demo{Proof}We first note that it will suffice to consider normalized
bases.  Suppose  that $(x_n)$ is  a normalized block basic sequence.
Let $a_n=\|x_n\|_2^2.$
Then for any finitely nonzero $(t_n)_{n\in \Cal M}$ with $\max|t_n|\le
1,$ we have:
$$
\align
 \Lambda (\sum_{n\in J}t_nx_n)&=
\sum_{n\in \Cal M}\sum_{j=1}^{\infty}F(|t_n||x_n(j)|)\\
&= \sum_{n\in\Cal
M}\sum_{j=1}^{\infty}(|t_n|^2F(|x_n(j)|)-|x_n(j)|^2|t_n|^2\log|t_n|)\\
&= \sum_{n\in\Cal M}|t_n|^2 (1- a_n\log |t_n|).
\endalign
$$
It follows easily that $(x_n)$ is 1-equivalent to the unit vector basis
in the Orlicz modular space $\ell(F_{a_n}).$

Now it is clear that if $a\le b\le Ma$ then $F_a(x)\le F_b(x)\le
MF_a(x)$ and from this it follows easily, using the uniform
$\Delta_2-$conditions on $F_t$ for $0\le t\le 1$ that there is constant
$K=K(M)$ so that if $a_n\le b_n\le Ma_n$ for $n\in\Cal M$ then the unit
vector bases of $\ell(F_{a_n})$ and $\ell(F_{b_n})$ are $K$-equivalent.

Now we turn to the general case.  First note that $\ell_F$ is
superreflexive and so if $X$ is any closed subspace and $2<p<\infty$ is
fixed
then $X^*$ is of cotype $p$ with some cotype constant $D$ independent of
$X.$   Now suppose $(x_n)_{n\in\Cal M}$ is a normalized $M$-unconditional
basic sequence in $\ell_F$ whose closed linear span is a subspace $X.$
Consider the co-ordinate functional $e_j^*$ as an element of $X^*$.  Then
$$ \sum_{n\in\Cal M}|x_n(j)|^p=\sum_{n\in\Cal M}|e_j^*(x_n)|^p \le
D^pM^p.$$

Now suppose $(t_n)_{n\in\Cal M}$ is finitely nonzero, and $\max|t_n|\le
1.$ Since
$\ell_F$ has cotype 2, there is (cf. \cite{15} Theorem 1.d.6) a universal
constant
$C$ so that
$$\frac1{CM} \|(\sum_{n\in\Cal M}|t_n|^2|x_n|^2)^{1/2}\|_F \le
\|\sum_{n\in\Cal M}t_nx_n\|_F\le
CM \|(\sum_{n\in\Cal M}|t_n|^2|x_n|^2)^{1/2}\|_F.$$

Let us calculate the modular $\Lambda(f)$ where $f=(\sum_{n\in\Cal
F}|t_n|^2|x_n|^2)^{1/2}.$  Then
$$ \Lambda(f)= \sum_{j=1}^{\infty}G(\sum_{n\in\Cal
M}|t_n|^2|x_n(j)|^2).$$
Since $G$ is concave, it is subadditive and so:
$$ \Lambda(f) \le \sum_{j=1}^{\infty}\sum_{n\in\Cal
M}G(|t_n|^2|x_n(j)|^2) = \sum_{n\in\Cal M}F_{a_n}(|t_n|).$$

For the reverse inequality consider
$$ \Lambda(M^{-1}D^{-1}f) =\sum_{j=1}^{\infty}G(\sum_{n\in\Cal
M}|t_n|^2w_n(j))$$
where $w_n(j)=M^{-2}D^{-2}|x_n(j)|^2.$  Then $(\sum_{n\in\Cal M
}w_n(j)^{p/2})^{2/p}\le 1.$  Thus we can apply Lemma 3.1 to deduce that
for each $j,$
$$ \sum_{n\in\Cal M}G(|t_n|^2w_n(j)) \le G(\sum_{n\in\Cal
M}|t_n|^2w_n(j))+ \frac{p}2\sum_{n\in\Cal M}G(w_n(j))|t_n|^2.$$

Summing over $j$, and using the fact that $MD\ge 1,$ we obtain
$$ \sum_{n\in\Cal M}F_{a_n}(M^{-1}D^{-1}|t_n|) \le \Lambda(f)
+ \sum_{n\in\Cal M}|t_n|^2.$$
The fact that $\ell_F$ has cotype 2 implies an estimate
$$ (\sum_{n\in\Cal M}|t_n|^2)^{1/2} \le CM\|\sum_{n\in\Cal M}t_nx_n\|_F.
$$

It follows easily that $(x_n)_{n\in\Cal M}$ is $K$-equivalent to the unit
vector basis of $\ell(F_{a_n})$ where $K$ depends only on $M.$  This and
the preceding remarks complete the proof.
 \enddemo

The following theorem follows immediately:

\proclaim{Theorem 4.3}Every unconditional basic sequence in
$\ell_F$ is equivalent to a sequence of constant coefficient
blocks in $\ell_F$ and hence spans a subspace isomorphic to a
complemented subspace of $\ell_F.$\endproclaim

Let us note that this implies a strong universality principle for
unconditional basic sequences in $\ell_F.$  Precisely, $\ell_F$ has
an unconditional basis (obtained by repeating every length constant
coefficient block infinitely often) so that every normalized
unconditional basic sequence is equivalent to a subsequence of the basis.
Such a property is also enjoyed by Pe\l czy\'nski's universal space
(\cite{14}, Theorem 2.d.10 or \cite{18}).
We next observe that $\ell_F$ obeys a strong form of the
Schroeder-Bernstein property for spaces with unconditional bases.

\proclaim{Theorem 4.4}Let $X$ be a Banach space with an unconditional
basis, and suppose that $X$ embeds into $\ell_F$ and $\ell_F$ embeds into
$X$.  Then $X$ is isomorphic to $\ell_F.$\endproclaim

\demo{Proof}By the preceding theorem $X$ is isomorphic to a complemented
subspace of $X$ spanned by constant coefficient blocks
$(u_n)_{n=1}^{\infty}.$  We now observe that $(u_n)_{n=1}^{\infty}$ must
contain an infinite number of blocks of the same length, for otherwise
$X$ is isomorphic to an Orlicz modular space $\ell_{F_{a_n}}$ where
$\lim_{n\to\infty}a_n=0$ and this can easily be seen not to contain a
copy of $\ell_F.$  Hence $\ell_F$ is complemented in $X$.  By Proposition
3.a.5 of \cite{14}, $\ell_F$ is isomorphic to $\ell_F\oplus X$ and this
is now trivially isomorphic to $X.$
\enddemo

In \cite{12} it is shown that any non-hereditarily Hilbertian space
with an unconditional basis contains a closed subspace with a 2-UFDD
which fails to have local unconditional strucure.  The following
theorem
(our main result of the
section) shows that this result cannot be substantially improved.

\proclaim{Theorem 4.5}Let $X$ be a closed subspace of $\ell_F$ with a
UFDD $(E_k)_{k=1}^{\infty}$ such that
the spaces $(E_k)$ are uniformly Hilbertian  (i.e. $\sup
d(E_k,\ell_2^{N_k})<\infty$, where $N_k=\dim E_k.$)  Then one can choose
an unconditional basis $(f_{ik})_{i=1}^{N_k}$ of $E_k$ so that the
collection
$(f_{ik})_{i,k}$ is an unconditional basis of $X.$\endproclaim

\demo{Remark}In particular the theorem applies to
any uniform-UFDD.\enddemo

\demo{Proof}Let $\|\,\|_{E_k}$ be a Euclidean norm on $E_k$ so that
$\|x\|_F\le \|x\|_{E_k}\le C\|x\|_F$ where $C=\sup d(E_k,\ell_2^{N_k}).$
Let $M$ be the constant of unconditionality for the Schauder
decomposition $(E_k).$
We choose a basis $(f_{ik})$ for $E_k$ which is orthonormal for both
$\|\,\|_{E_k}$ and $\|\,\|_2.$  Suppose $(t_{ik})$ is finitely non zero
and that $(\epsilon_{ik})$ is a choice of signs.
Let $x_k =\sum_{i=1}^{N_k}t_{ik}f_{ik}$ and
$y_k=\sum_{i=1}^{N_k}\epsilon_{ik}t_{ik}f_{ik}.$  Then for the set
$\Cal M$ of all
$k$ such that they are nonzero we have $\|x_k\|_F/\|y_k\|_F\le C$ and
$\|y_k\|_F/\|x_k\|_F\le C.$  We also have $\|x_k\|_2=\|y_k\|_2$.  Both
$(x_k)_{k\in\Cal M}$ and $(y_k)_{k\in \Cal M}$ are $M$-unconditional
basic sequences.  Hence they are $K$-equivalent by Lemma 4.2 where
$K=K(C,M).$  In particular,
 $$ \|\sum_{k\in\Cal M}y_k\|_F \le K\|\sum_{k\in\Cal M}x_k\|_F$$
whence the basis $(f_{ik})$ is $K$-unconditional.\enddemo

\demo{Remark}The properties of unconditional basic sequences in $\ell_F$
have other applications, for example to uniqueness questions.  We plan to
discuss these applications in a separate paper.
\enddemo

\enddocument